\documentclass[11pt]{amsart}
\usepackage{amssymb,amsmath,txfonts,mathrsfs,mathtools,titletoc, tikz, stmaryrd}
\usepackage{color}
\usepackage{graphicx}
\usepackage{float}
\usepackage{appendix}
\usepackage{hyperref}

\usepackage[alphabetic]{amsrefs}
\newtheorem{theorem}{Theorem}[section]

\newtheorem{lemma}[theorem]{Lemma}
\newtheorem{remark}[theorem]{Remark}

\newtheorem{definition}[theorem]{Definition}
\newtheorem{cor}[theorem]{Corollary}

\newtheorem{conj}[theorem]{Conjecture}

\numberwithin{equation}{section}

\def\pf{{\it Proof:}~}

\usepackage{color}
\usepackage{graphicx}
\usepackage{float}
\usepackage{appendix}

\usepackage[alphabetic]{amsrefs}
\numberwithin{equation}{section}

\usepackage[alphabetic]{amsrefs}

\numberwithin{equation}{section}

\def\pf{{\it Proof:}~}

\def\pf{{\it Proof:}~}

\numberwithin{equation}{section}

\def\pf{{\it Proof:}~}

\begin{document}

\title[The sharp diameter bound of diameter]{The sharp diameter bound of stable minimal surfaces}
\author{Qixuan Hu, Guoyi Xu, Shuai Zhang}
\address{Qixuan Hu\\ Department of Mathematical Sciences\\Tsinghua University, Beijing\\P. R. China}
\email{huqx20@mails.tsinghua.edu.cn}
\address{Guoyi Xu\\ Department of Mathematical Sciences\\Tsinghua University, Beijing\\P. R. China}
\email{guoyixu@tsinghua.edu.cn}
\address{Shuai Zhang\\ Department of Mathematical Sciences\\Tsinghua University, Beijing\\P. R. China}
\email{zhangshu22@mails.tsinghua.edu.cn}
\date{\today}
\date{\today}

\begin{abstract}
For three dimensional complete Riemannian manifolds with scalar curvature no less than one, we prove that Schoen-Yau diameter estiamte for stable minimal surfaces in those manifolds is sharp. 
\\[3mm]
\end{abstract}
\thanks{G. Xu was partially supported by NSFC 12141103.}

\maketitle

\section{Introduction}

Let $(M^3, g)$ be a complete $3$-dim Riemannian manifold with scalar curvature $R\geq 1$, for any stable minimal surface $\Sigma$ in $M$ (possibly $\partial \Sigma\neq \emptyset$), and $\gamma\subset \Sigma$ be a closed curve, if 
\begin{align}
U_\rho(\gamma)\cap \partial\Sigma= \emptyset, \quad \quad \mathrm{Image}[H_1(\gamma)\rightarrow H_1(U_\rho(\gamma))]\nequiv 0. 
\end{align}
where $U_\rho(\gamma)\vcentcolon= \{x\in \Sigma: \mathrm{dist}_{\Sigma}(x, \gamma)\leq \rho\}$, Gromov and Lawson \cite[Theorem $10.2$]{GL} showed $\rho\leq\pi$. Furthermore they \cite[Remark $10.6$]{GL} conjecture that the best possible conclusion would be $\rho\leq\frac{\pi}{\sqrt{2}}$ which is achieved by $M=\mathbb{S}^2(\sqrt{2})\times \mathbb{S}^1$, $\Sigma=\mathbb{S}^2(\sqrt{2})\subset M$, and $\gamma$ be the great circle in $\Sigma$.

On the other hand, Schoen-Yau \cite[Proof of Proposition $1$]{SY-CMP} showed that $\rho\leq \frac{\sqrt{6}}{3}\pi$; which implies that the diameter of any complete stable minimal surface in a complete Riemannian manifold $(M^3, g)$ with scalar curvature $R(g)\geq 1$ has upper bound $\frac{2\sqrt{6}\pi}{3}$.

In this note, we show Schoen-Yau's diameter bound for stable minimal surfaces is sharp, which is a corollary of the following theorem.

\begin{theorem}\label{thm diam of stable main result}
{Let $\Sigma$ be a complete stable minimal surface in a complete Riemannian manifold $(M^3, g)$ with scalar curvature $R(g)\geq 1$. Then 
\begin{align}
\mathrm{Diam}(\Sigma)< \frac{2\sqrt{6}\pi}{3}. \label{sharp ineq of diam of stable main}
\end{align}

The upper bound is sharp in the following sense: there is a sequence of complete $3$-dim manifolds $M_k\vcentcolon= (S^2\times S^1, g_k)$ with $R(g_k)\geq 1$ and $k\in \mathbb{Z}^+$, and compact stable minimal surfaces $\Sigma_k\subseteq M_k$ such that $\displaystyle\lim_{k\rightarrow\infty}\mathrm{Diam}(\Sigma_k)= \frac{2\sqrt{6}}{3}\pi$.
}
\end{theorem}

The strictness of (\ref{sharp ineq of diam of stable main}) is obtained by a simple observation, the main contribution of this note is the construction of the example manifolds, which shows the sharpness of this upper bound.

From \cite[4.C]{Gro20}, we know
\begin{align*}
\mathrm{FillRad}(M)\leq\frac{1}{2}\mathrm{width}_{n-1}(M)\leq \frac{1}{2}\mathrm{Diam}(M)
\end{align*}
for complete Riemannian manifold $(M^n,g)$, where $\mathrm{FillRad}(M)$ is the filling radius and $\mathrm{width}_{n-1}(M)$ is the $(n-1)$-th Uryson width, defined as in \cite[section 4]{Gro20}. 

The results in section \ref{sec 3dim case} show that
for a complete stable minimal surface $\Sigma$ in a complete $3$-dim Riemannian manifold $(M^3, g)$ with scalar curvature $R\geq 1$, we have $\displaystyle \mathrm{FillRad}(\Sigma)< \frac{\sqrt{6}\pi}{3}$. 

Gromov \cite{Gro07} proposed the following conjecture on the filling radius.
\begin{conj}
    If $(M^n, g)$ is a complete Riemannian manifold with scalar curvature $R(g)\geq \sigma^2>0$, then there exists a universal constant $c_n$ depending only on $n$ such that \begin{align*}
        \mathrm{FillRad}(M)\leq \frac{c_n}{\sigma}.
    \end{align*}
\end{conj}
The conjecture has been partially answered by \cite{WXYZ} for manifolds with finite asymptotic dimension. 

In $3$-dim complete manifolds with scalar curvature no less than $-1$, Munteanu, Sung, Wang \cite{MSW} obtained the area upper bound of stable minimal surfaces, 

In \cite[3.10]{Gro23} Gromov conjectured that a complete manifold with
scalar curvature $R(g) \geq 6$ admits a singular foliation by surfaces of area and diameter bounded by a universal constant. The compact case was solved in \cite{LM}, and the non-compact case was solved in \cite{LW}; also see \cite{WZ} when $M$ has boundary. 

The organization of this paper is as follows. In section \ref{sec surface case}, we firstly get a general sharp diameter upper bound, for surfaces with the positive first  eigenvalue corresponding to Laplace operator with suitable curvature potential term. The example Riemannian surfaces are also constructed in this section. The key is to get the suitable function from some special ODE (see (\ref{new ODE for psi})), which comes from the symmetrization of the corresponding PDE along the diameter (or `radius') direction. 

In section \ref{sec 3dim case}, we use the functions from the example Riemannian surfaces to construct the corresponding $3$-dim complete Riemannian manifolds containing the stable minimal surfaces, whose diameter approximates the sharp upper bound.

\section{The first eigenvalue and the diameter}\label{sec surface case}

\begin{definition}\label{def the first eigenvalue of complete mfld}
{For a complete Riemannian manifold $(M^n, g)$ and $w\in C^\infty(M^n)$, we define
\begin{align}
\lambda_1(-\Delta+ w)\vcentcolon= \inf_{f\in H_c^{1}(M^n)}\frac{\int_M |\nabla f|^2+ w\cdot f^2}{\int_M f^2}, \nonumber 
\end{align}
where $H_c^{1}(M)$ is the closure of $C_c^\infty(M^n)$ in $H^{1}(M^n)$. 
}
\end{definition}

\begin{lemma}\label{lem posi eigenvalue implies positive eigenfunc}
{For a complete Riemannian manifold $(M^n, g)$ with $w\in C^\infty(M^n), \lambda\in \mathbb{R}^+$, assume $\lambda_1(-\Delta+ w)\geq \lambda$, then there is $\varphi\in C^\infty(M^n)$ satisfying
\begin{align}
\varphi> 0, \quad \quad \text{and} \quad \quad -\Delta \varphi+ (w- \lambda)\cdot \varphi\geq 0. \nonumber 
\end{align}
}
\end{lemma}

\pf
{By definition of $\lambda_1$, we get 
\begin{align}
\lambda_1(-\Delta+(w-\lambda))&= \inf_{f\in H_c^{1, 2}(M^n)}\frac{\int_M |\nabla f|^2+ (w-\lambda)\cdot f^2}{\int_M f^2}= \inf_{f\in H_c^{1, 2}(M^n)}\frac{\int_M |\nabla f|^2+ w\cdot f^2}{\int_M f^2}-\lambda\nonumber\\
&=\lambda_1(-\Delta+w)-\lambda\geq 0.
\end{align}

The compact case follows from similar argument of \cite[Theorem $8.38$]{GT-PDE}. The non compact case follows from \cite[Theorem $1$]{FS}.
}
\qed

\begin{theorem}\label{thm sharp upper bound of diameter}
{For any $\beta> \frac{1}{4}$ and $\lambda> 0$, assume $\lambda_1(-\Delta_\Sigma+ \beta\cdot K_\Sigma)\geq \lambda$ on a complete $2$-dim Riemannian manifold $(\Sigma, g)$, where $K_\Sigma$ is the sectional curvature of $(\Sigma, g)$. Then 
\begin{align}
\mathrm{Diam}(\Sigma)< \frac{2\beta \pi}{\sqrt{\lambda\cdot (4\beta- 1)}}. \nonumber 
\end{align}

The upper bound is sharp in the following sense: there is a sequence of complete Riemannian manifolds $\Sigma_k= (S^2, g_k)$ with $\lambda_1(-\Delta_{\Sigma_k}+ \beta\cdot K_{\Sigma_k})\geq \lambda$ and
\begin{align}
\lim_{k\rightarrow \infty}\mathrm{Diam}(\Sigma_k)= \frac{2\beta \pi}{\sqrt{\lambda\cdot (4\beta- 1)}}. \nonumber 
\end{align}
}
\end{theorem}

\begin{remark}\label{rem PDE Bonnet-Myers thm}
{Most part of the following proof is similar to \cite[Lemma $2.4$]{WXZ} (also see \cite[Appendix A]{XuKai}) except the last part, we include the complete proof here for readers' convenience. For $\beta= 1$, Theorem \ref{thm sharp upper bound of diameter} also answers a question raised in \cite[Remark $2.5$]{WXZ}.  
}
\end{remark}

\pf
{\textbf{Step (1)}. By Lemma \ref{lem posi eigenvalue implies positive eigenfunc}, we get there exists a positive function $v\in C^\infty(M)$, such that 
\begin{align}
-\Delta v+(\beta\cdot K_\Sigma-\lambda)v\geq0. \nonumber 
\end{align}

Let $u=v^{\frac1\beta}$, we get
\begin{align}
-\Delta u+(K_\Sigma-\lambda\beta^{-1})u+(1-\beta)u^{-1}|\nabla u|^2\geq0.\label{beta version eq of u}
\end{align}

Case(i):If $\Sigma$ is compact, let $p,q\in \Sigma$ be two points with largest distance. Define $I[\gamma]$ by
\begin{align}
I[\gamma]=\int_\gamma u \mathrm{d}s,\nonumber 
\end{align}
where $\mathrm{d}s$ is arclength along $\gamma$ and $\gamma$ is any curve from $p$ to $q$.

Case(ii): If $\Sigma$ is non compact, chose $p\in\Sigma$ and $R= \frac{2\beta \pi}{\sqrt{\lambda\cdot (4\beta- 1)}}$, Define $I[\gamma]$ by
\begin{align}
I[\gamma]=\int_\gamma u \mathrm{d}s,\nonumber 
\end{align}
where $ds$ is arclength along $\gamma$ and $\gamma$ is any curve from $p$ to $\partial B_R(p)$.

Since $u$ is positive, for any curve $\gamma:[0,L]\to M$ connecting $p$ and $\partial B_R(p)$, if there exists $L_1\in(0,L)$, such that $\gamma(L_1)\in\partial B_R(p)$. Let $\gamma_1$ be the restriction of $\gamma$ on $[0,L_1]$, we have
\begin{align}
I[\gamma]=\int_\gamma u ds>\int_{\gamma_1}u ds=I[\gamma_1]. \nonumber 
\end{align}

So
\begin{align}
\inf_{\gamma\subset M}I[\gamma]=\inf_{\mathrm{Int}(\gamma)\subset B_R(p)}I[\gamma], \nonumber 
\end{align}
where $\mathrm{Int}(\gamma)$ is the interior of $\gamma$.

In both cases, there are at least one minimizer of $I[\cdot]$. Let $\gamma_0:[0,l]\to M$ be one of the minimizers of $I[\cdot]$ and parametrized with unit speed.

\textbf{Step (2)}.  Let $V(s)= \varphi(s)\cdot \vec{n}$ be the variation vector field along $\gamma_0(s)$, where $\vec{n}$ is the unit normal vector field along $\gamma$ and 
\begin{align}
\varphi(s)= u^{-\frac{1}{2}}(\gamma_0(s))\psi(s), \quad \quad \psi(s)= \sin\Big(\frac{\pi}{l}s\Big). \label{choice of varphi and psi}
\end{align}

The first variation formula yields
\begin{align}
\langle \nabla_{\partial s}\partial s, \vec{n}\rangle= \frac{\partial}{\partial\vec{n}}\ln u.\label{mean curvature vanish}
\end{align}

The non-negativity of the second variation of $I$ at $\gamma_0$ gives
\begin{align}
0\leq \delta^2I[\gamma_0]= \int_{\gamma_0} [D^2u(\vec{n}, \vec{n})- 2u\cdot (\frac{\partial}{\partial\vec{n}}\ln u)^2]\cdot \varphi^2+ u\cdot \Big((\varphi')^2- K_\Sigma\cdot \varphi^2\Big)\mathrm{d}s. \label{2nd vari}
\end{align}

Note $\displaystyle \Delta_\Sigma u= D^2u(\vec{n}, \vec{n})+ D^2u(\partial s, \partial s)$, combining (\ref{mean curvature vanish}), we get
\begin{align}
\Delta_\Sigma u&= D^2u(\vec{n}, \vec{n})+ D^2u(\partial s, \partial s)= D^2u(\vec{n}, \vec{n})+ u''- (\frac{\partial}{\partial\vec{n}}\ln u)^2\cdot u, \label{Hess of u}
\end{align}
where  $u'\vcentcolon= \frac{\mathrm{d}}{\mathrm{d} s}u(\gamma_0(s))$.

Plugging \eqref{Hess of u} into \eqref{2nd vari}, using \eqref{beta version eq of u}, we obtain
\begin{align}
0&\leq \int_{\gamma_0} (\Delta_\Sigma u- u''- u\cdot (\frac{\partial}{\partial\vec{n}}\ln u)^2)\cdot \varphi^2+ u\cdot \Big((\varphi')^2- K_\Sigma\cdot \varphi^2\Big)\mathrm{d}s \nonumber  \\
&\leq \int_{\gamma_0} \Big((1-\beta)|\nabla u|^2u^{-1}- u''- u\cdot (\frac{\partial}{\partial\vec{n}}\ln u)^2\Big)\cdot \varphi^2+ u\cdot \Big((\varphi')^2- \lambda\beta^{-1} \varphi^2\Big)\mathrm{d}s. \nonumber \\
&\leq \int_{\gamma_0} \Big((1-\beta)(u')^2u^{-1}- u''\Big)\cdot \varphi^2+ u\cdot \Big((\varphi')^2- \lambda\beta^{-1} \varphi^2\Big)\mathrm{d}s.\label{1-dim case}
\end{align}

Plugging \eqref{choice of varphi and psi} into \eqref{1-dim case}, and use integration by part, we get
\begin{align}
0&\leq\int_{\gamma_0}\Big[(1-\beta)(u')^2u^{-2}\psi^2+u'(-u^{-2}u'\psi^2+2u^{-1}\psi\psi')+u(-\frac12u^{-\frac32}u'\psi+u^{-\frac12}\psi')^2-\lambda\beta^{-1}\psi^2\Big]\mathrm{d}s \nonumber \\
&=\int_{\gamma_0}\Big[(1+\frac14(\beta-\frac14)^{-1})(\psi')^2-\lambda\beta^{-1}\psi^2-\Big((\beta-\frac14)^{\frac12}u'u\psi-\frac12(\beta-\frac14)^{-\frac12}\psi'\Big)^2\Big]\mathrm{d}s\nonumber\\
&\leq\int_{\gamma_0}\Big[(1+\frac14(\beta-\frac14)^{-1})(\psi')^2-\lambda\beta^{-1}\psi^2\Big]\mathrm{d}s\nonumber\\
&=\int_0^l\Big[(1+\frac14(\beta-\frac14)^{-1})(\frac{\pi}{l})^2\cos^2(\frac{\pi}{l}s)-\lambda\beta^{-1}\sin^2(\frac{\pi}{l}s)\Big]\mathrm{d}s\nonumber\\
&=\Big[\frac{4\beta}{4\beta-1}(\frac{\pi}{l})^2-\lambda\beta^{-1}\Big]\cdot\frac{l}{2}.\nonumber
\end{align}

So we get $l\leq\frac{2\beta\pi}{\sqrt{\lambda\cdot(4\beta-1)}}$. If $l=\frac{2\beta\pi}{\sqrt{\lambda\cdot(4\beta-1)}}$, we get 
\begin{align}
\int_0^l\Big((\beta-\frac14)^{\frac12}u'u\psi-\frac12(\beta-\frac14)^{-\frac12}\psi'\Big)^2\mathrm{d}s=0.
\end{align}
Note $(\beta-\frac14)^{\frac12}u'u\psi-\frac12(\beta-\frac14)^{-\frac12}\psi'$ is a continuous function on $[0,l]$, so 
\begin{align}
(\beta-\frac14)^{\frac12}u'u\psi-\frac12(\beta-\frac14)^{-\frac12}\psi'\equiv0,\quad \forall s\in[0,l].
\end{align}
However, for $s=0$, we have $\psi(0)=0$ and $\psi'(0)=\frac{\pi}{l}$, which is the contradiction. 

Now we have $l< \frac{2\beta\pi}{\sqrt{\lambda\cdot(4\beta-1)}}$. If $\Sigma$ is not compact, from the definition of $l$, we know that
\begin{align}
l\geq d(p, \partial B_R(p))= \frac{2\beta\pi}{\sqrt{\lambda\cdot(4\beta-1)}}, \nonumber 
\end{align}
which is the contradiction. 

Hence $\Sigma$ must be compact, and from the definition of $l$ in compact case we get 
\begin{align}
\mathrm{Diam}(\Sigma)< \frac{2\beta\pi}{\sqrt{\lambda\cdot(4\beta-1)}}. \nonumber 
\end{align}

\textbf{Step (3)}. In the rest argument, we construct $\Sigma_k$.  We assume $k\in \mathbb{Z}^+$ with 
\begin{align}
k^{-1}\in \Big(0, \frac{\pi}{4}\sqrt{\frac{\beta}{\lambda}}\Big). \label{choice of epsilon}
\end{align}

Define 
\begin{align}\nonumber 
\phi(x)=\left\{\begin{array}{lll}
e^{-\frac{1}{x}},&\quad &x>0,\\
0,&\quad &x\leq0.
\end{array}\right.
\end{align}
\begin{align}
\eta_k(x)&= \frac{\phi(2k^{-1}-x)}{\phi(x-k^{-1})+\phi(2k^{-1}-x)},\quad \quad \forall x\in \mathbb{R}. \label{def of eta-1}
\end{align}
Note $\eta_k$ is a smooth function on $\mathbb{R}$.

We define
\begin{equation}\label{new ODE for psi}
\left\{
\begin{array}{rl}
&\psi_k'(t)= -\frac{1}{\beta}\psi_k^2+ \frac{1- \eta_k^2}{4\beta^2}\psi_k^2- \lambda, \quad \quad \forall t\in [k^{-1}, b_k), \\
&\psi_k(t)=\sqrt{\beta\lambda}\cdot\cot(\sqrt{\frac{\lambda}{\beta}}t),\quad \quad \forall t\in (0, k^{-1}],
\end{array} \right.
\end{equation} 
where $[k^{-1}, b_k)$ is the largest interval for the solution to (\ref{new ODE for psi}). 

By (\ref{new ODE for psi}), note $\beta> \frac{1}{4}$, we get 
\begin{align}
\psi_k'(t)< 0, \quad \quad \quad \forall t\in [0, b_k). \label{psi is decreasing}
\end{align}

Using (\ref{new ODE for psi}) , note $|\eta_k|\leq 1$, we get
\begin{equation}\label{ODE psi-initial}
\left\{
\begin{array}{rl}
&\psi_k'(t)\geq -\frac{1}{\beta}\psi_k^2- \lambda \ ,  \quad \quad \forall  t\in[k^{-1},2k^{-1}], \\
&\psi_k(k^{-1})= \sqrt{\beta\lambda}\cdot\cot(\sqrt{\frac{\lambda}{\beta}}k^{-1}). 
\end{array} \right.
\end{equation} 

From (\ref{ODE psi-initial}), taking integral from $k^{-1}$ to $2k^{-1}$ with respect to $t$, we get
\begin{align}
\psi_k(2k^{-1})\geq \sqrt{\lambda\beta}\cot\Big\{\sqrt{\frac{\lambda}{\beta}}k^{-1}+ \cot^{-1}(\frac{\psi_k(k^{-1})}{\sqrt{\lambda\beta}})\Big\}=\sqrt{\lambda\beta}\cot\Big(2\sqrt{\frac{\lambda}{\beta}}k^{-1}\Big)> 0, \label{lower bound of psi 2epsilon}
\end{align}
where the last inequality follows from (\ref{choice of epsilon}). 

\textbf{Step (4)}. Now we consider $\psi_k(t)$ for $t\in [2k^{-1}, b_k)$, and we have
\begin{align}
\psi_k'(t)= -(\frac{1}{\beta}- \frac{1}{4\beta^2})\psi_k^2- \lambda \ ,  \quad \quad \forall  t\in [2k^{-1}, b_k).  \label{the last part ODE for psi}
\end{align}

Note $\beta> \frac{1}{4}$, hence $\frac{1}{\beta}- \frac{1}{4\beta^2}> 0$. Direct integration yields
\begin{align}
\psi_k(t)= 2\beta\sqrt{\frac{\lambda}{4\beta- 1}}\cdot \cot\Big\{\frac{\sqrt{(4\beta- 1)\lambda}}{2\beta}(t- 2k^{-1})+ \cot^{-1}(\frac{1}{2\beta}\sqrt{\frac{4\beta- 1}{\lambda}}\psi_k(2k^{-1}))\Big\}, \quad  \forall t\in [2k^{-1}, b_k). \label{keq eq of psi t}
\end{align}

From (\ref{keq eq of psi t}), we get $b_k< \infty$ and $\displaystyle \lim_{t\rightarrow b_k-}\psi_k(t)= -\infty$. 

Define $l_k= \min\{t\in (0, b_k): \psi_k(t)= 0\}$, from (\ref{lower bound of psi 2epsilon}) and (\ref{psi is decreasing}), we get 
\begin{align}
l_k\geq 2k^{-1}, \quad \quad \quad \psi_k(l_k)= 0. \nonumber 
\end{align}

From (\ref{keq eq of psi t}) we get
\begin{align}
\frac{\sqrt{(4\beta- 1)\lambda}}{2\beta}(l_k- 2k^{-1})+ \cot^{-1}(\frac{1}{2\beta}\sqrt{\frac{4\beta- 1}{\lambda}}\psi_k(2k^{-1}))= \frac{\pi}{2}, \nonumber 
\end{align}
which implies
\begin{align}
l_k&= \frac{2\beta}{\sqrt{(4\beta- 1)\lambda}}\Big\{\frac{\pi}{2}- \cot^{-1}(\frac{1}{2\beta}\sqrt{\frac{4\beta- 1}{\lambda}}\psi_k(2k^{-1}))\Big\}+ 2k^{-1} \nonumber \\
&= \frac{\beta \pi}{\sqrt{(4\beta- 1)\lambda}}- \frac{2\beta}{\sqrt{(4\beta- 1)\lambda}}\cot^{-1}(\frac{1}{2\beta}\sqrt{\frac{4\beta- 1}{\lambda}}\psi_k(2k^{-1}))+ 2k^{-1}. \label{eq of l}
\end{align}


From (\ref{lower bound of psi 2epsilon}) and (\ref{eq of l}), we have
\begin{align}
&l_k\geq \frac{\beta \pi}{\sqrt{(4\beta- 1)\lambda}} - \frac{2\beta}{\sqrt{(4\beta- 1)\lambda}}\cot^{-1}\Big(\sqrt{1- \frac{1}{4\beta}}\cot\Big(2\sqrt{\frac{\lambda}{\beta}}k^{-1}\Big)\Big)+ 2k^{-1}. \label{lower bound of l}
\end{align}

One key fact is that
\begin{align}
\varliminf_{k\rightarrow \infty}l_k\geq \frac{\beta \pi}{\sqrt{(4\beta- 1)\lambda}}. \label{limit of l wrt epsilon}
\end{align}

\textbf{Step (5)}. Now we define $f_k: [0, 2l_k]\rightarrow \mathbb{R}$ by
\begin{equation}\label{definition of f} 
\left\{
\begin{array}{lll}
(\ln f_k)'(t)= \frac{2\beta+ \eta_k(t)- 1}{2\beta^2}\psi_k(t) \ ,&  \quad  &\forall t\in (0, l_k];\\
f_k(k^{-1})=\sqrt{\frac{\beta}{\lambda}}\sin(\sqrt{\frac{\lambda}{\beta}}k^{-1})\ ;& & \\
f_k(t)=f_k(2l_k-t)\ ,& &\forall t\in[l_k,2l_k); \\
f_k(0)= f_k(2l_k)= 0. & &
\end{array} \right.
\end{equation}
From the property of $f_k$, we get the complete Riemannian surface $\Sigma_k= (S^2, dr^2+ f_k(r)^2d\theta)$ for $r\in [0, 2l_k]$, where $\theta\in \mathbb{S}^1$.

Now we define $\zeta_k: [0, 2l_k]\rightarrow \mathbb{R}$ by
\begin{equation}\label{definition of zeta}
\left\{
\begin{array}{lll}
(\ln \zeta_k)'(t)= \frac{1- \eta_k(t)}{2\beta}\psi_k(t) \ ,  &\quad & \forall t\in (0, l_k];\\
\zeta_k(k^{-1})=1 \ ; & &\\
\zeta_k(t)=\zeta_k(2l_k-t)\ , & & \forall t\in[l_k,2l_k); \\
\zeta_k(0)= \zeta_k(2l_k)=1.  & &
\end{array} \right.
\end{equation}
Note $\psi_k(t)\geq 0$ for $t\in (0, l_k]$ by the choice of $l_k$, then from (\ref{definition of zeta}), we get $\zeta_k(t)> 0$ for any $t\in [0, 2l_k]$.

Define $u_k(r, \theta)= \zeta_k(r)$, from (\ref{definition of zeta}) we know that $u_k\in C^\infty(\Sigma_k)$. From (\ref{new ODE for psi}), (\ref{definition of f}) and (\ref{definition of zeta}), we get that for $r\in(0,l_k]$, 
\begin{align}
&-\Delta_{\Sigma_k} u_k+ \beta\cdot K_{\Sigma_k} u_k- \lambda u_k= -(\zeta_k''+ \frac{f_k'}{f_k}\zeta_k')+ \beta(-\frac{f_k''}{f_k})\zeta_k- \lambda\zeta_k \label{eq of u} \\
=&-\zeta_k\Big\{\psi_k'+\frac{1}{\beta}\psi_k^2-\frac{1- \eta_k^2}{4\beta^2}\psi_k^2+\lambda\Big\}= 0. \nonumber 
\end{align}

By the symmetry definition of $u_k,f_k$, we get \eqref{eq of u} also holds for $r\in[l_k,2l_k)$. So
\begin{align}
-\Delta_{\Sigma_k} u_k+\beta\cdot K_{\Sigma_k} u_k=\lambda u_k. \nonumber 
\end{align} 

From \eqref{limit of l wrt epsilon} and Theorem \ref{thm sharp upper bound of diameter}, we have 
\begin{align}
\frac{2\beta \pi}{\sqrt{(4\beta- 1)\lambda}}\geq \varlimsup_{k\rightarrow \infty}\mathrm{Diam}(\Sigma_k)\geq \varliminf_{k\rightarrow \infty}\mathrm{Diam}(\Sigma_k)\geq 2\varliminf_{k\rightarrow \infty}l_k\geq \frac{2\beta \pi}{\sqrt{(4\beta- 1)\lambda}}. \nonumber 
\end{align}
}
\qed

\begin{remark}\label{rem shape of example mflds}
{We shall point out that in our example, for the case of $\beta=\frac12$, from (\ref{new ODE for psi}) we have
\begin{align*}
    \psi_k(t)\leq \psi_k(k^{-1})=\sqrt{\frac{\lambda}{2}}\cdot\cot(\sqrt{2\lambda}k^{-1})\leq \frac{k}{2},\quad t\in(k^{-1},2k^{-1}].
\end{align*}
From (\ref{definition of f}) we have

\begin{equation}
\left\{
\begin{array}{lll}
(\ln f_k)'(t)\leq 2\psi_k(t)\leq k \ ,&  \quad  &\forall t\in (k^{-1}, 2k^{-1}],\\
(\ln f_k)'(t)= 0 \ ,&  \quad  &\forall t\in (2k^{-1}, l_k],\\
f_k(k^{-1})=\sqrt{\frac{1}{2\lambda}}\sin(\sqrt{2\lambda}k^{-1}).& & \\
\end{array} \right.
\end{equation}
So
\begin{align*}
    f_k(t)\equiv f_k(2k^{-1})\leq e\cdot f_k(k^{-1})\leq e\cdot k^{-1},\quad t\in [2k^{-1},2l_k-2k^{-1}].
\end{align*}
and $\Sigma_k=(S^2,g_k)$ converges to a segment $[0, 2\pi]$ as $k\to\infty$.
So $f_k$ converges to $0$ as $k\to\infty$, which implies $\Sigma_k$ converges to a segment (see Figure \ref{fig exmaple for dim 2}).

\begin{figure}[H]
    \centering
    \includegraphics[width=0.4\linewidth]{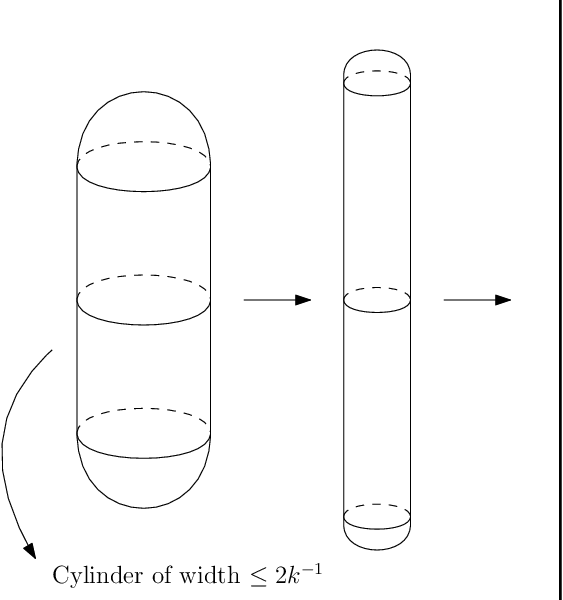}
    \caption{Figure of $\Sigma_k$ as $k\to\infty$.}
    \label{fig exmaple for dim 2}
\end{figure}
}
\end{remark}

As an corollary, we get a sharp version of \cite[Lemma $16$]{CL}. 
\begin{cor}\label{cor the diameter upper bound of 2-dim surface-CL}
{For compact surface $(\Sigma^2, g)$, if there is a smooth function $u>0$ such that 
\begin{align}
\frac{\Delta u}{u}\leq (K_\Sigma- 2^{-1})+\frac{1}{2}\frac{|\nabla_\Sigma u|^2}{u^2}, \label{diff ineq of u} 
\end{align} 
then the diameter of $\Sigma$ satisfies the sharp upper bound $\mathrm{Diam}(\Sigma)< 2\pi$.
}
\end{cor}

\pf
{Let $v= \sqrt{u}> 0$, from (\ref{diff ineq of u}), we get
\begin{align}
-\Delta v+ (\frac{1}{2}K_\Sigma- \frac{1}{4})\cdot v\geq 0. \label{v diff ineq} 
\end{align}

Using (\ref{v diff ineq}) in the proof of Theorem \ref{thm sharp upper bound of diameter}, we get the conclusion. 
}
\qed

\section{The diameter of stable minimal surfaces in $3$-dim manifolds}\label{sec 3dim case}

\begin{lemma}\label{lem stable ms implies diff ineq of u}
Let $\Sigma$ be a complete stable minimal surface in a complete $3$-dim Riemannian manifold $(M^3, g)$ with scalar curvature $R\geq 1$, then $\displaystyle \lambda_1(-\Delta_\Sigma+ K_\Sigma)\geq \frac{1}{2}$, where $K_\Sigma$ is the sectional curvature of $\Sigma$. 

\end{lemma}

\pf
{Let $e_1,e_2,e_3$ be an orthonormal frame defined locally on $\Sigma$ with $e_1,e_2$ tangential and $e_3$ be unit normal. Since $\Sigma$ is a stable minimal surface, for all $f\in C_c^\infty(\Sigma)$ we have 
\begin{align}\label{stable equ}
\int_\Sigma|\nabla f|^2-(Rc(e_3,e_3)+\sum h_{ij}^2)f^2\geq0,
\end{align}
where $h_{ij}=<\nabla_{e_i}e_j,e_3>$ is the second fundamental form of $\Sigma$. (see \cite[Chapter $1$]{LiBook})

Since $\Sigma$ is minimal, we have $h_{11}+h_{22}=0$. By Gauss curvature equation, we have $K_\Sigma=R_{1212}+h_{11}h_{22}-h_{12}^2$, and $R_{ijkl}$ is the Ricci curvature tensor of $M$. 

So \eqref{stable equ} can be written as
\begin{equation}
\int_{\Sigma}|\nabla f|^2-(\frac12R-K_\Sigma+\frac12\sum h_{ij}^2)f^2\geq 0. \nonumber 
\end{equation}

Since $R\geq1$, we get
\begin{equation}
\int_{\Sigma}|\nabla f|^2+(K_\Sigma-\frac12)f^2\geq 0, \nonumber 
\end{equation}
which implies $\lambda_1(-\Delta+ K_\Sigma)\geq \frac{1}{2}$. 
}
\qed

\begin{theorem}\label{thm diam of stable ms}
{Let $\Sigma$ be a complete stable minimal surface in a complete $3$-dim Riemannian manifold $(M^3, g)$ with scalar curvature $R\geq 1$, then
\begin{align}
\mathrm{Diam}(\Sigma)< \frac{2\sqrt{6}\pi}{3}. \label{sharp ineq of diam of stable ms}
\end{align}

The upper bound is sharp in the following sense: there is a sequence of complete $3$-dim manifolds $M_k\vcentcolon= (S^2\times S^1, g_k)$ with $R\geq 1$ and $k\in \mathbb{Z}^+$, and compact stable minimal surfaces $\Sigma_k\subseteq M_k$ such that $\displaystyle\lim_{k\rightarrow\infty}\mathrm{Diam}(\Sigma_k)= \frac{2\sqrt{6}}{3}\pi$.
}
\end{theorem}

\pf
{\textbf{Step (1)}. The inequality \eqref{sharp ineq of diam of stable ms} follows from Lemma \ref{lem stable ms implies diff ineq of u} and Theorem \ref{thm sharp upper bound of diameter} for the case $\beta=1, \lambda= \frac{1}{2}$.

Let $g_k= \mathrm{d}r^2+f_k^2(r)\mathrm{d}\theta^2+\zeta_k^2(r)\mathrm{d}\phi^2$, where $(r,\theta)$ is the polar coordinate on $S^2$, $\phi$ is the coordinate on $S^1$, and $f_k$, $\zeta_k$ are defined as in \eqref{definition of f} and \eqref{definition of zeta}. 

Note $r\in [0, 2l_k]$, then by (\ref{limit of l wrt epsilon}) and (\ref{sharp ineq of diam of stable ms}), we have 
\begin{align}
\lim_{k\rightarrow\infty}\mathrm{Diam}(\Sigma_k)= \frac{2\sqrt{6}}{3}\pi. \nonumber 
\end{align}

Define $\Sigma_k=S^2\times\{\phi_0\}\subseteq M_k$ in the rest argument. 

If the $r$-coordinate of $p\in M_k$ is $0$ or $2l_k$, by the definition of $M_k$, the $k^{-1}$-neighborhood of $p$ in $M_k$ is an open set in $\mathbb{S}^2(\sqrt{2})\times \mathbb{S}^1(1)$; where $\mathbb{S}^2(\sqrt{2})$ is a round sphere with radius $\sqrt{2}$, $\mathbb{S}^1(1)$ is a round unit circle and $\mathbb{S}^2(\sqrt{2})\times \mathbb{S}^1(1)$ has the product metric. 

The $k^{-1}$-neighborhood of $p$ in $\Sigma$ is an open set of $\mathbb{S}^2(\sqrt{2})\times\{\phi_0\}\subset \mathbb{S}^2(\sqrt{2})\times \mathbb{S}^1(1)$ for some $\phi_0\in \mathbb{S}^1$. 

Since $\mathbb{S}^2(\sqrt{2})$ is a minimal surface of $\mathbb{S}^2(\sqrt{2})\times \mathbb{S}^1(1)$ with the product metric, we get the mean curvature of $\Sigma_k$ is $0$ in a neighborhood of $p$.

\textbf{Step (2)}. Now we assume $p\in \Sigma_k$, where $\phi_0\in S^1$, and the $r$-coordinate of $p$ is not equal to $0$ or $2l_k$.

Let $\{e_1,e_2,e_3\}=\{\frac{\partial}{\partial r},\frac{1}{f_k}\frac{\partial}{\partial\theta},\frac{1}{\zeta_k}\frac{\partial}{\partial\phi}\}$ be a local orthonormal frame on a neighborhood of $p$ in $M_k$, we have $e_1,e_2$ is tangential and $e_3$ is unit normal. The second fundamental form of $\Sigma_k$ is defined by symmetric quadratic tensor $h_{ij}=<\nabla_{e_i}e_j,e_3>$, $i,j\in\{1,2\}$, where $\nabla$ is the Riemannian connection of $M_k$. 

Let $\Gamma_{ij}^k=\frac12g^{lk}\left(\frac{\partial g_{jl}}{\partial x_i}+\frac{\partial g_{il}}{\partial x_j}-\frac{\partial g_{ij}}{\partial x_l}\right)$ be the Christoffel symbol, where $\frac{\partial}{\partial x_1},\frac{\partial}{\partial x_2},\frac{\partial}{\partial x_3}$ are $\frac{\partial}{\partial r},\frac{\partial}{\partial \theta},\frac{\partial}{\partial\phi}$ respectively, and $(g_{ij})$ is the matrix corresponding to the metric of $M_k$.

By direct computation, we get $\Gamma_{33}^1=-\zeta_k\frac{\partial \zeta_k}{\partial r}$, $\Gamma_{13}^3=\Gamma_{31}^3=\frac{1}{\zeta_k}\frac{\partial \zeta_k}{\partial r}$, $\Gamma_{22}^1=-f_k\frac{\partial f_k}{\partial r}$, and $\Gamma_{12}^2=\Gamma_{21}^2=\frac{1}{f_k}\frac{\partial f_k}{\partial r}$. Other Christoffel symbols are zero.

We have 
\begin{align}
\nabla_{e_1}e_2&=\nabla_{\frac{\partial}{\partial r}}(\frac{1}{f_k}\frac{\partial}{\partial\theta})=\frac{\partial}{\partial r}(\frac{1}{f_k})\frac{\partial}{\partial\theta}+\frac{1}{f_k^2}\frac{\partial f_k}{\partial r}\frac{\partial}{\partial\theta}=0, \nonumber \\
\nabla_{e_1}e_1&=\nabla_{\frac{\partial}{\partial r}}\frac{\partial}{\partial r}=0, \nonumber \\
\nabla_{e_2}e_2&=\nabla_{(\frac{1}{f_k}\frac{\partial}{\partial\theta})}(\frac{1}{f_k}\frac{\partial}{\partial\theta})
=\frac{1}{f_k}\frac{\partial}{\partial\theta}(\frac{1}{f_k})\frac{\partial}{\partial\theta}-\frac{1}{f_k}\frac{\partial f_k}{\partial r}\frac{\partial}{\partial r}. \nonumber 
\end{align}

We get $\nabla_{e_i}e_j$ is in the tangent plane of $\Sigma_k$, and 
\begin{align}
h_{ij}=0, \quad \quad \quad \forall i, j\in \{1,2\}. \nonumber 
\end{align}
So on a neighborhood of $p$ in $\Sigma_k$, we have the mean curvature of $\Sigma_k$ is $0$.

So we get $\Sigma_k$ is a minimal surface.

\textbf{Step (3)}. Since $h_{ij}=0$, we get the scalar curvature of $M_k$ is $2K_{\Sigma_k}+2Rc(e_3,e_3)$, where $K_{\Sigma_k}$ is the sectional curvature of $\Sigma_k$. 

By the definition of Ricci curvature tensor, we have 
\begin{align}
R_{ij}= \frac{\partial}{\partial x_k}\Gamma_{ij}^k- \frac{\partial}{\partial x_j}\Gamma_{ik}^k+ \Gamma_{ij}^s\Gamma_{sk}^k- \Gamma_{ik}^s\Gamma_{sj}^k. \nonumber 
\end{align}

We get 
\begin{align}
Rc(e_3,e_3)=\frac{1}{\zeta_k^2}Rc(\frac{\partial}{\partial \phi},\frac{\partial}{\partial \phi})=-\frac{1}{\zeta_k}\cdot(\frac{\partial^2\zeta_k}{\partial r^2}+\frac{1}{f_k}\frac{\partial f_k}{\partial r}\frac{\partial \zeta_k}{\partial r}). \nonumber 
\end{align}

On $\Sigma_k$, we have $\Delta_{\Sigma_k}=\frac{\partial^2}{\partial r^2}+\frac{f_k'}{f_k}\frac{\partial}{\partial r}+\frac{1}{f_k^2}\frac{\partial^2}{\partial\theta^2}$, and $K_{\Sigma_k}=-\frac{f_k''}{f_k}$, where $\Delta_{\Sigma_k}$ is the Laplacian operator on $\Sigma_k$. 

Define $u_k(r, \theta)= \zeta_k(r)$, then $u_k\in C^\infty(\Sigma_k)$. So $Rc(e_3,e_3)=-\frac{\Delta_{\Sigma_k} u_k}{u_k}$. 

By \eqref{eq of u}, we get
\begin{align}
-\Delta_{\Sigma_k} u_k+(K_{\Sigma_k}-\frac{1}{2})u_k=0.\label{eq of zeta}
\end{align}

We get the scalar curvature of $M_k$ is 
\begin{align}
R(g_k)=2K_{\Sigma_k}-2\frac{\Delta_{\Sigma_k}u_k}{u_k}=1. \nonumber 
\end{align}

\textbf{Step (4)}. Let $w_k=\ln u_k$, from \eqref{eq of zeta} we get
\begin{align}
\Delta_{\Sigma_k} w_k=\frac{\Delta_{\Sigma_k} g}{g}-|\frac{\nabla g}{g}|^2=K_{\Sigma_k}-\frac12-|\nabla w_k|^2. \label{def of w}
\end{align}

For any $\varphi\in C^\infty(\Sigma_k)$, we multiply both sides of (\ref{def of w}) by $\varphi^2$ and take integration on $\Sigma_k$ to get 
\begin{align}
\int_{\Sigma_k} -2\varphi\nabla \varphi\cdot\nabla w_k=\int_{\Sigma_k}(K_{\Sigma_k}-\frac12)\varphi^2-\varphi^2\cdot|\nabla w_k|^2. \label{inte ineq with varphi and w}
\end{align}

Note 
\begin{align}
-2\varphi\nabla \varphi\cdot\nabla w_k\geq -|\nabla\varphi|^2-\varphi^2\cdot|\nabla w_k|^2. \label{CS ineq}
\end{align}

By (\ref{CS ineq}) and (\ref{inte ineq with varphi and w}), we obtain
\begin{align}
\int_{\Sigma_k}|\nabla \varphi|^2+(K_{\Sigma_k}-\frac{1}{2})\varphi^2\geq 0,\quad \quad  \forall \varphi\in C^\infty(\Sigma_k). \nonumber 
\end{align}

Since $K_{\Sigma_k}-\frac12=\frac{\Delta_{\Sigma_k}u_k}{u_k}=-Rc(e_3,e_3)$ and $h_{ij}=0$, we get 
\begin{align}
\int_{\Sigma_k}|\nabla \varphi|^2-(Rc(e_3,e_3)+\sum h_{ij}^2)\varphi^2\geq 0, \quad \quad \forall \varphi\in C^\infty(\Sigma_k). \nonumber 
\end{align}

So $\Sigma_k$ is a stable minimal surface of $M_k$.
}
\qed

\section*{Acknowledgments}
We thank Bo Zhu for his comments and suggestion on the earlier version of this paper.  

\textbf{Data availability} Data sharing not applicable to this article as no datasets were generated
or analysed during the current study.

\textbf{Declarations}

\textbf{Conflict of interests} The authors declare that they have no conflict of interest.

\end{document}